\documentclass[12pt]{article}
\usepackage{amsmath,amsfonts,amssymb,amscd}
\usepackage{graphicx}
\usepackage[T2A]{fontenc}

\newtheorem{theo}{Theorem}
\newtheorem{lem}{Lemma}[section]

\newtheorem{cor}{Corollary}[section]

\makeatletter \@addtoreset{equation}{section} \makeatother

\newcommand{\mR}{\mathbb{R}}
\newcommand{\mT}{\mathbb{T}}
\newcommand{\mZ}{\mathbb{Z}}

\newcommand{\one}{{\bf 1}}

\newcommand{\bw}{{\bf w}}

\newcommand{\calA}{{\cal A}}

\newcommand{\calT}{{\cal T}}
\newcommand{\calU}{{\cal U}}

\newcommand{\thet}{\vartheta}

\newcommand{\id}{\operatorname{id}}

\newcommand\qed{{\unskip\nobreak\hfil\penalty50
  \hskip2em\hbox{}\nobreak\hfil\mbox{\rule{1ex}{1ex} \qquad}
    \parfillskip=0pt \finalhyphendemerits=0\par\medskip}}

\begin{document}

\title
{Volume preserving diffeomorphisms as Poincare maps for volume preserving flows}
\author{D.Treschev \\
Steklov Mathematical Institute of Russian Academy of Sciences
}
\date{}
\maketitle

\begin{abstract}
Let $q:M\to M$ be a volume-preserving diffeomorphism of a smooth manifold $M$. We study the possibility to present $q$ as the Poincar\'e map, corresponding to a volume-preserving vector field on $\mT\times M$, $\mT = \mR/\mZ$.
\end{abstract}


\section{Introduction}

Dynamical systems are split essentially in two classes: systems with continuous time and the ones with discrete time. These classes are closely related and there are canonical tools to transmit results from one class to another. The passage from continuous time to discrete one is performed through the Poincar\'e map while the inverse passage is usually associated with various forms of the suspension construction.

These two tool are not at all straightforward and various technical problems appear if one tries to implement them. One of the main problems is as follows. The conventional suspension construction, \cite{BS} associates with a diffeomorphism $q: M\to M$ of a smooth manifold $M$ a flow on a certain manifold $\widehat M$,
$\dim\widehat M = \dim M + 1$, whose topology is a priori not clear. Applications usually require $\widehat M\cong \mT\times M$, where $\mT = \mR/\mZ$ is a circle, hence the suspension should be constructed in another way.

The situation is especially clear and simple if instead of the total class of systems with continuous time we restrict to the subclass of non-autonomous systems with periodic time dependence. In this context most of results on the existence of suspensions with $\widehat M = \mT\times M$ are known, see for example, \cite{TZ}, Section 1.3 and references therein.

In this note we present some results on the existence of a suspension on $\mT\times M$
in the case when $q$ is a volume preserving diffeomorphism. Our motivation was an attempt to answer questions of authors of \cite{KKP} who need the corresponding results in the study of mixing properties for motions of an ideal fluid.

\section{Main result}

All the objects below will be $C^\infty$-smooth.

Let $M$ be a compact manifold, $\dim M = m$, and let $\nu$ be a volume form on $M$. Here we mean that $\nu$ is a differential $m$-form nowhere vanishing and such that
$\int_M \nu >0.$ We denote $\mT=\mR/\mZ=\{t\bmod 1\}$. Then
\begin{equation}
\label{omega}
  \omega = dt\wedge\nu
\end{equation}
is a volume form on the manifold $\mT\times M$.

Let $\pi_\mT:\mT\times M\to\mT$ and let $\pi_M:\mT\times M\to M$ be the natural projections
$$
  \mT\times M\ni (t,x) \mapsto \pi_\mT(t,x) = t, \quad
  (t,x) \mapsto \pi_M(t,x) = x.
$$
Consider the vector field $v$ on $\mT\times M$. We assume that
\medskip

(A) the first component of $v$ is positive: $D\pi_\mT\, v = v_\mT > 0$,

(B) $v$ preserves the form $\omega$: $L_v\omega=0$, where $L_v$ is the Lie derivative.
\medskip

Let $g_v^s$ be the flow, generated by the vector field $v$ on $\mT\times M$. Condition (A) implies that the Poincar\'e map (the first return map)
$P_v:\{0\}\times M\to\{0\}\times M$ is well-defined:
$$
  \{0\}\times M \ni x \mapsto P_v(x) = g_v^{\sigma_v(x)}(x),
$$
where $\sigma_v : \{0\}\times M \to \mR_+ = \{s\in\mR : s>0\}$ is defined by
$$
  \sigma_v(x) = \min_{s>0} \{s : \pi_\mT g_v^s(x) = 0 \}.
$$
The map $P_v$ preserves the volume form
$\lambda = i_v\omega |_{\{t=0\}}$ on $\{0\}\times M$.

\begin{theo}
\label{theo:inclu}
Let $Q:\{0\}\times M\to \{0\}\times M$ be another map which preserves $\lambda$. We assume that $Q$ is (smoothly) isotopic to $P_v$ in the group of $\lambda$-preserving self maps of $\{0\}\times M$. Then there exists an $\omega$-preserving vector field $u$ on $\mT\times M$ such that $D\pi_\mT\, u > 0$ and $Q = P_u$.
\end{theo}

\section{Preliminary remarks}

Consider the smooth in $s$ family of diffeomorphisms $\calT_M(s):M\to M$ of a manifold $M$ into itself, $s\in\mR$. We extend $\calT_M$ to a family of diffeomorphisms $\calT^s$ of $\mR\times M$ putting by definition
\begin{equation}
\label{calT}
      \mR\times M
  \ni (t,x) \mapsto \calT^s(t,x)
   =  \Big(t+s , \calT_M(t+s)\circ\calT_M^{-1}(t)(x)\Big).
\end{equation}
Direct computation shows that $\calT^s$ is a flow i.e.,
$$
  \calT^0 = \id \quad\mbox{and}\quad
  \calT^{s_2} \circ \calT^{s_1} = \calT^{s_1+s_2}\quad
  \mbox{ for any $s_1,s_2\in\mR$}.
$$
The flow $\calT^s$ generates the vector field $\calA$ on $\mR\times M$:
\begin{equation}
\label{calA}
  \calA = \Big( \frac d{ds}\Big|_{s=0} \calT^s \Big) \circ \calT^{-s} , \qquad
  D\pi_\mR(\calA) = 1,
\end{equation}
where $\pi_\mR : \mR\times M \to \mR$ is a natural projection.

Let $S:\mR\times M\to \mR\times M$ be the time-one shift:
$$
  (t,x) \mapsto S(t,x) = (t+1,x).
$$
The flow $\calT^s$ is said to be periodic (with period 1) if it commutes with $S$.

If $\calT^s$ is periodic then the field $\calA$ is also periodic. Indeed,
\begin{eqnarray*}
     \calA\circ S
 &=& \Big(\frac d{ds}\calT^s\Big)\circ\calT^{-s}\circ S
  =  \frac d{ds}\Big(S\circ\calT^s\circ S^{-1}\Big)\circ\calT^{-s}\circ S \\
 &=& DS\Big(\frac d{ds}\calT^s\Big)\circ S^{-1}\circ\calT^{-s}\circ S
  =  DS\,\calA.
\end{eqnarray*}
In this case the canonical projection $\mR\times M\to (\mR\times M)/S = \mT\times M$ generates a flow $T^s$ and a vector field $\alpha$ on $\mT\times M$, the projections of $\calT^s$ and $\calA$ to $\mT\times M$.
\medskip

Now we discuss the property of volume preservation.

\begin{lem}
\label{lem:1}
Suppose the vector field $\alpha$ on the manifold $N$ preserves the volume form $\mu$. Then

1. $d\imath_\alpha \mu = 0$,

2. for any positive function $f:N\to\mR$ the vector field $\displaystyle\frac1f \alpha$ preserves the form $f\mu$.
\end{lem}

{\it Proof of Lemma \ref{lem:1}}. The first statement follows from the Cartan formula $L_\alpha = d\imath + \imath d$. The second statement follows from the computation
$$
     L_{\frac1f \alpha} f\mu
  =  d\imath_{\frac1f \alpha} f\mu
  =  d\imath_\alpha \mu
  =  L_\alpha \mu
  = 0.
$$
\qed

\begin{lem}
\label{lem:form}
Let $\calT_M(s) : M\to M$ be a smooth in $s\in\mR$ family of diffeomorphisms, let $\mu_s$ be a smooth in $s$ family of volume forms on $M$ such that
$$
  \calT_M^*(s) \mu_s = \mu_0,
$$
and let $\kappa_t = \pi_M^*\mu_t|_{M_t}$ be the volume form on $M_t = \{ t\}\times M$ (the natural push up of $\mu_t$).

Then the flow $\calT^s$ preserves the form $dt\wedge\kappa_t$.
\end{lem}

{\it Proof}. It is sufficient to note that $(\calT^s)^*\kappa_t = \kappa_t$ and  $(\calT^s)^* dt = dt + \beta$, where the 1-form $\beta$ does not contain $dt$. \qed

\section{Proof of Theorem \ref{theo:inclu}}
\label{sec:proof}

(a) We put $\displaystyle\widehat v = \frac1{v_\mT} v$. Then $D\pi_\mT\,\widehat v = 1$. Let
$$
  (t,x) \mapsto g_{\widehat v}^s(t,x), \qquad
  (t,x)\in \mT\times M
$$
be the flow of the vector field $\widehat v$. By Lemma \ref{lem:1} $g_{\widehat v}^s$ preserves the form $v_\mT \omega$. The Poincar\'e maps $P_v$ and
$P_{\widehat v}$ coincide. Hence,
\begin{equation}
\label{ph=P}
  g_{\widehat v}^1(0,x) = (1,P_v(x))\qquad
  \mbox{for any $x\in M$}.
\end{equation}

Below we need a lift of the flow $g_{\widehat v}^s$ to the covering space $\mR\times M$ to a periodic flow $G^s$. We determine the family $\sigma_s : M\to M$ by
$$
  G^s(0,x) = (s,\sigma_s(x)), \qquad
  s\in\mR.
$$

(b) Let $\gamma_s$ be a smooth in $s$ isotopy from conditions of Theorem \ref{theo:inclu}: for any $s\in [0,1]$ the map $\gamma_s$ is a $\lambda$-preserving diffeomorphism of $M \cong \{0\}\times M$,
\begin{equation}
\label{gamma(01)}
  \gamma_0 = P_v , \quad  \gamma_1 = Q.
\end{equation}
Changing smoothly parametrization on $\gamma_s$, we can assume that $\gamma_s = P_v$ in a neighborhood of $\{s=0\}$ and $\gamma_s = Q$ in a neighborhood of $\{s=1\}$. We extend $\gamma_s$ to all the axis $\mR = \{s\}$, for example, putting $\gamma_s = P_v$ for $s<0$ and  $\gamma_s = Q$ for $s>0$.

(c) Consider the family of maps $\calT_M(s):M\to M$,
$$
    \calT_M(s)
  = \sigma_s\circ P_{\widehat v}^{-1}\circ \gamma_s, \qquad
     s\in\mR.
$$
By (\ref{ph=P}) and (\ref{gamma(01)})
\begin{equation}
\label{thet_preserves}
  \calT_M(0) = \id, \quad
  \calT_M(1) = Q, \qquad
  \mbox{hence $\calT_M(1)$ preserves the form $\lambda$}.
\end{equation}

Let $\calT^s$ be the flow on $\mR\times M$ generated by the family $\calT_M(s)$ and let $\calU$ be the corresponding vector field on $\mR\times M$. Then by (\ref{calT})
$$
  D\pi_\mT\,\calU = 1, \quad
       \calT^0 = \id_{\mT\times M}, \quad
  \calT^1(0,x) = (1,Q(x)).
$$

By (\ref{calT}) and (\ref{calA})
$$
    \calU
  = \bigg( 1, \frac d{ds}\Big|_{s=0}
                \Big( \sigma_{t+s}\circ P_{\widehat v}^{-1}\circ\gamma_{t+s} \Big)
           \circ \Big(\sigma_t\circ P_{\widehat v}^{-1}\circ\gamma_t \Big)^{-1}
    \bigg)
  = \widehat v + \bw,
$$
where
\begin{equation}
\label{vw}
    \widehat v
  = \bigg( 1, \Big( \frac d{dt} \sigma_t \Big)
              \circ \sigma_t^{-1} \bigg),   \quad
    \bw
  = \bigg( 0, D(\sigma_t\circ P_{\widehat v}^{-1})
              \Big( \frac d{dt} \gamma_t \Big)
           \circ \gamma_t^{-1}\circ P_{\widehat v}\circ\sigma_t^{-1} \bigg) .
\end{equation}

Near the points $t=0$ and $t=1$ we have: $d\gamma_t / dt = 0$. Therefore the vector field $\calU$ coincides with $\widehat v$. Hence $\calU|_{s\in[0,1]}$ can be extended to a periodic vector field $\widehat\calU$ on $\mR\times M$. Let $\widehat\thet^s$ be the corresponding flow on $\mR\times M$. Due to periodicity of $\widehat\thet^s$ and $\widehat\calU$ their projections to the flow $\thet^s$ and the vector field $U$ on $\mT\times M$ are well-defined.

(d) Let $\one$ be the vector field on $\mR\times M$ determined by the equations
$$
  D\pi_\mT\one = 1, \quad    D\pi_M\one = 0.
$$
By Lemma \ref{lem:form} the flow $\widehat\thet^s$ preserves some volume form $\Omega$ on $\mR\times M$ which can be chosen so that
\begin{equation}
\label{iOmega}
  \imath_\one\Omega|_{\{t=0\}} = \lambda.
\end{equation}

Any volume form on $\mT\times M$ equals $\widehat\rho\omega$. where $\widehat\rho:M\to\mR$ is a function. Therefore $\Omega=\widehat\rho\omega$, where by Lemma \ref{lem:form}, (\ref{thet_preserves}), and (\ref{iOmega})
$$
    \imath_\one\Omega|_{t=0}
  = \imath_\one\Omega|_{t=1}
  = \lambda.
$$
This equation and periodicity of $\widehat\thet^s$ imply that the function $\widehat\rho$ is 1-periodic in $t$. Hence there exists a function
$\rho:\mT\times M\to\mR$ such that $\widehat\rho = \rho\circ\pi_\mT$. The flow $\thet^s$ preserves the volume form $\rho\omega$.

By Lemma \ref{lem:1} the vector field $u = \rho U$ preserves $\omega$. It remains to note that $P_u = P_U = Q$. \qed

\section{Generalizations and remarks}

{\bf 1. The case when $Q$ is a small perturbation of $P_v$}

\begin{cor}
By (\ref{vw}) the vector field $\bw$ is small if the isotopy $\gamma_s$ is near-identity. Moreover, in this case $U$ is close to $\widehat v$, therefore $\rho$ is close to $v_\mT$ and $u$ is a small perturbation of $v$.
\end{cor}

Here the following question appears. Suppose two $\lambda$-preserving diffeomorphisms $Q,P_v : M\to M$ are close to each other. Is it true that there exists a $\lambda$-preserving isotopy $\gamma_s$, joining $Q$ with $P_v$, such that for any $s\in[0,1]$ the diffeomorphism $\gamma_s$ lies in a small neighborhood of $P_v$?

In general we do not know the answer. However in the case $\dim M = 2$ the form $\lambda$ may be identified with a symplectic structure, and the standard argument, using a proper generating function, gives a positive answer. Details can be found, for example, in \cite{ST}, proof of Lemma 1.
\medskip

{\bf 2. The case of a finite smoothness}

Suppose that the form $\nu$ and the vector field $v$ are only $C^k$-smooth. Suppose also that the isotopy $\gamma_s$, joining $P_v$ with $Q$, consists of $C^k$-smooth diffeomorphisms and depends $C^{k+1}$-smoothly on $s$. Then straightforward translation of the argument from Section \ref{sec:proof} to the spaces $C^k$ gives a result analogous to Theorem \ref{theo:inclu} in finitely smooth categories: more precisely, the vector field $u$, which is constructed in Theorem \ref{theo:inclu}, belongs to the space $C^k$.

\medskip

{\bf 3. The real-analytic case}

If the manifold $M$, the vector field $v$, and the diffeomorphism $Q$ are real-analytic, then assuming that the isotopy $\gamma_s$ is a curve in the group of real-analytic $\lambda$-preserving maps, we can construct $u$ also real-analytic. Here the dependence of $\gamma_s$ on $s$ need not be real-analytic: it can be only $C^1$-smooth. Indeed, by using the continuous averaging procedure from \cite{PT} (see also \cite{TZ}), we can transform $\gamma_s$ to another isotopy between $P_v$ and $Q$ which is real-analytic and $1$-periodic in $s$.
\medskip

{\bf 4. What happens if $M$ is not closed?}

Suppose $M$ equals the closure $\overline D$ of an open domain $D$, lying in a (may be, non-compact) manifold $N$. We assume that

(a) $M$ is compact,

(b) $P_v(M) = M$,

(c) There exists an open domain $D^+\supset M$ such that the maps $Q : D^+\to N$ and
$P_v|_{D^+} : D^+\to N$ are smoothly isotopic to each other in the class of $\lambda$-preserving maps from $D^+$ to $N$.

Then the same argument imply the existence of an $\omega$-preserving vector field $u$ on $\mT\times N$ such that $D\pi_\mT\, u > 0$ and $Q|_{\overline D} = P_u|_{\overline D}$.
Moreover, if $Q$ = $P_v$ in a neighborhood of the boundary $\partial D$ of
$M =\overline D$, then $u$ equals $v$ in this neighborhood.

\end{document}